\documentclass[12pt,a4paper]{article}
\usepackage[top=3cm, bottom=3cm, left=3cm, right=3cm]{geometry}
\usepackage[latin1]{inputenc}
\usepackage{amsmath}
\usepackage{amsthm}
\usepackage{amsfonts}
\usepackage{amssymb}
\usepackage{graphics}
\usepackage[hang,flushmargin]{footmisc} 

\usepackage{amsmath,amsthm,amssymb,graphicx, multicol, array}
\usepackage{enumerate}
\usepackage{enumitem}
\usepackage{xcolor}
\usepackage{currfile}

\usepackage[bookmarks,colorlinks,breaklinks, pagebackref]{hyperref}

\hypersetup{linkcolor=blue,citecolor=red,filecolor=blue,urlcolor=blue} 

\usepackage{tabto}

\numberwithin{equation}{section}

\usepackage{tikz}
\usepackage{pgfplots}

\DeclareMathOperator{\li}{li}

\newtheorem{thm}{Theorem}[section]
\newtheorem{lem}{Lemma}[section]
\newtheorem{conj}{Conjecture}[section]

\newtheorem{cor}{Corollary}[section]

\usepackage{fancyhdr}
\title{Note On Prime Gaps And Zero Spacings}
\date{}
\author{N. A. Carella}

\begin{document}
\maketitle

\textbf{\textit{Abstract}:} 
The note proves several related problems for the greatest lower bound liminf, and the least upper bound
limsup of the zero spacings. \let\thefootnote\relax\footnote{\today \date{} 
}

\tableofcontents

\section{Introduction } \label{sec1}
The paper focuses on the problems of prime gaps and zero spacings. It draws from the extensive literature on
the theory of prime gaps and zero spacings and a prime and zero duality principle to achieve possible solutions
of several related problems such as the greatest lower bound $\liminf$, the least upper bound $\limsup$ of the zero
spacings, and the least upper bound limsup of the prime gaps (known as Cramer conjecture).\\

An apparently simple new idea linking the prime numbers and the zeros of the zeta function, see Section \ref{sec5}, and
recent results in the literature of analytic number theory will be utilized to achieve the followings results.\\

\begin{thm} \label{thm1.1} There exists infinitely many pairs of consecutive zeros $\rho_n=\beta_n+i\gamma_n$ and $\rho_{n+1}=\beta_{n+1}+i\gamma_{n+1}$ of the
	zeta function $\zeta(s)$ on the critical line $\Re e(s)=1/2$ such that
	\begin{equation}
		\liminf_{n\to \infty} \frac{\gamma_{n+1}-\gamma_n}{2 \pi /\log(\gamma_n/2 \pi)}=0.
	\end{equation}
	and 
	\begin{equation}
		\limsup_{n\to \infty} \frac{\gamma_{n+1}-\gamma_n}{2 \pi /\log(\gamma_n/2 \pi)}=\infty.
	\end{equation}
\end{thm}

The same analysis provides the concrete bounds stated below.\\

\begin{cor}    \label{cor1.1} Let $\{1/2+\gamma_n : n\geq 1 \}$ be the subset of zeros of $\zeta(s)$ along the critical line $\Re e(s) = 1/2$. Then\\
	{\normalfont(i)} There is an infinite subsequence of small zero spacings such that
	\begin{equation}
		\gamma_{n+1}-\gamma_n\leq c_0\frac{(\log \log \gamma_n)^2}{(\log \gamma_n)^{3/2}}.
	\end{equation}
	
	{\normalfont(ii)} There is an infinite subsequence of large zero spacings such that
	\begin{equation}
		\gamma_{n+1}-\gamma_n\geq c_1\frac{(\log \log \gamma_n)(\log \log \log \log \gamma_n)}{(\log \gamma_n) (\log \log \log \gamma_n)^2},
	\end{equation}
	where $c_0, c_1 > 0$ are constants.
\end{cor}

\begin{cor}   \label{cor1.2} For a large integer $n \geq 1 $, the following statements hold.\\ 
	
	{\normalfont(i)} Unconditionaly, the $n$th prime gap satisfies the asymptotic upper bound 
	\begin{equation}
		2 \leq p_{n+1}-p_n = O\left (\log^2 p_n \right ).
	\end{equation}
	{\normalfont(ii)} Assume the Riemann hypothesis. Then, 
	\begin{equation}
		2 \leq p_{n+1}-p_n = O\left (\frac{\log^2 p_n}{\log \log n} \right ).
	\end{equation}
\end{cor}

\begin{proof}[\textbf{Proof}] (i) Consider the strictly positive quantity $p_{n+1}-p_n\geq 2$. Use Lemma \ref{lem5.1} to obtain
	\begin{eqnarray} 
		2 &\leq & p_{n+1}- p_n  \\
		&=& \left (   \gamma_{n+k}-\gamma_n  \right )\frac{ \log ^2 n}{2 \pi} \left (1+O\left ( \frac{\log \log n}{ \log n}\right )\right ) \nonumber \\
		&=& O \left ( \log^2 n \right ) \nonumber .
	\end{eqnarray} This unconditional upper bound follows from $\gamma_{n+k}-\gamma_n \ll 1$. (ii) The conditional upper bound follows from Theorem \ref{thm4.2}.
\end{proof}

\textbf{Note 2.} All these results assume that every nontrivial zero of the zeta
function $\zeta(s)$ is of the form $\rho= 1/2 + i\gamma$. However, since
\begin{eqnarray}\label{100}
	 |\rho_{n+1}-\rho_n|&=&|\beta_{n+1}-\beta_n+i(\gamma_{n+1}-\gamma_n)| \\
	 &\leq& 1+\gamma_{n+1}-\gamma_n\nonumber\\
&\ll& 1\nonumber,	 
\end{eqnarray}
these results are valid even if some nontrivial zeros are not on the critical line $\Re e(s) = 1/2$.\\

Section \ref{sec2} provides an introduction to various conjectures in the literature and establishes the notation used
throughout the paper. The technical background required in the new results are given in Section \ref{sec3} to Section \ref{sec5}. The proof of
Theorem \ref{thm1.1} appears in Section \ref{sec6}.

\section{Various Conjectures} \label{sec2}
The $n$th zero of the zeta function $\zeta(s)$ in the critical strip $0 < \Re e(s) < 1$ is denoted by $\rho_n = \beta_n + i\gamma_n$ with $\beta_n, \gamma_n \in \mathbb{R}$, and the $n$th prime in the sequence of primes $
2, 3, 5, \ldots,$ is denoted by $p_n$.\\

\subsection{Average Spacing And Normalized Spacing} Let $\{ \rho_n = \beta_n + i\gamma_n : n\geq 1\}$ be the subsequence of zeros of the
zeta function $\zeta(s)$ on the critical strip $0 <\Re e(s) < 1$. The subset $\{ 1/2 +\gamma_n : n\geq 1 \}$ of zeros of $\zeta(s)$ coincides
with the sequence of real zeros of the entire function $Z(t)=e^{i \theta(t)} \zeta(1/2+i t)$, where $\theta(t)= \arg \zeta(1/ 2+ it)$. The imaginary $\gamma_n$ parts are assumed to form an increasing sequence
$14 \leq \gamma_1 \leq \gamma_2 \leq \cdots \leq \gamma_n \leq \cdots $ of real numbers. The average zero spacing between the critical zeros $\rho=\eta + it$ of $\zeta(s)$ in the rectangle $\{ s = \beta+ it : 0\leq \beta \leq 1; 0\leq  t
\leq T \}$ is the ratio of the maximum height $T > 0$ to the number $N(T)$ of zeros inside the aforementioned rectangle.
Specifically, the average zero spacing has the aymptotic formula
\begin{equation}
	\frac{T}{N(T)}=\frac{2 \pi}{\log(T / 2 \pi)}\left (1-\frac{1}{\log(T / 2 \pi)}+O\left ( \frac{1}{T}\right ) \right ),
\end{equation}
see Theorem \ref{thm4.1}. The normalized zero spacing is defined by
\begin{equation}
	\delta_n=(\gamma_{n+1}-\gamma_n)\frac{\log(T / 2 \pi)}{ 2 \pi}.
\end{equation}

The sequence of normalized zero spacings has the average value given by the sum \begin{equation}\label{105}
\overline{\delta_x}=\sum_{n \leq x}(\gamma_{n+1}-\gamma_n)\frac{\log(T / 2 \pi)}{ 2 \pi}=1. 
\end{equation}
For more details on a more general framework
for the sequence of zeros of entire functions and the corresponding probability densities, refer to \cite{LA2006}.
The normalized spacings of the nontrivial zeros are believed to be random as specified by the eigenvalues
distribution of $N\times N$ random Hermitian matrices for sufficiently large $N$. This arrangement matches $k$
consecutive zeros to $k$ eigenvalues of an $N\times N$ random Hermitian matrix, see \cite{MT1973}, \cite{OD1987} and similar papers for
further details.\\

The limit infimum and limit supremum of the sequence of normalized zero spacings $\{ \delta_n= \gamma_{n+1}-\gamma_n : n\geq 1 \}$ of
the zeros of $\zeta(s)$ are defined by
\begin{equation}
	\mu= \liminf_{n\to \infty} \frac{\gamma_{n+1}-\gamma_n}{2 \pi /\log(\gamma_n/2 \pi)} \quad 
	\text{     and     }
	\quad \lambda=\limsup_{n\to \infty} \frac{\gamma_{n+1}-\gamma_n}{2 \pi /\log(\gamma_n/2 \pi)},
\end{equation}
respectively. Here several conjectures are listed to guide the readers on the current state of knowledge on these
topics. A few proven results are also listed.\\

\begin{conj}  {\normalfont (\cite{MT1973})}  \label{conj2.1}  The following limits hold:
	\begin{equation}
		\liminf_{n\to \infty} \frac{\gamma_{n+1}-\gamma_n}{2 \pi /\log(\gamma_n/2 \pi)}=0 \quad 
		\text{\normalfont     and     }
		\quad \limsup_{n\to \infty} \frac{\gamma_{n+1}-\gamma_n}{2 \pi /\log(\gamma_n/2 \pi)} =\infty\nonumber.
	\end{equation}
\end{conj}
These are also known as the \textit{Small Gap Conjecture} and \textit{Large Gap Conjecture} respectively.\\

Assuming the random matrices conjectures and the Riemann hypothesis, the limit supremum has been
proved in \cite{ST2008}. But the only unconditional result in the literature on the existence of large spaces between the
zeros of the zeta function is the following.\\

\begin{thm}  {\normalfont (\cite{HL2008})}   \label{thm2.1}  Let $\rho_n = \beta_n + it_n$ be the $n$th critical zero of $\zeta(s)$ such that $0 \leq t_i < t_j$ for $i < j$. Then, the limit supremum of the normalized spacing satisfies the lower bound
	\begin{equation}
		\limsup_{n\to \infty} \frac{\gamma_{n+1}-\gamma_n}{2 \pi /\log(\gamma_n/2 \pi)}\geq 2.630637\ldots \nonumber .
	\end{equation}
\end{thm}

On the other end of the spectrum, there are a few results on the least upper bound for liminf similar to the
followings one.\\

\begin{thm}  {\normalfont (\cite{CG1984})}    \label{thm2.2} Assume the Riemann hypothesis. Then 
	\begin{equation}\liminf_{n\to \infty} \frac{\gamma_{n+1}-\gamma_n}{2 \pi /\log(\gamma_n/2 \pi)} \leq .5172\ldots \nonumber .
	\end{equation} 
\end{thm}

There are several equivalent definitions of the concept of pair correlation of the zeros of a function and its
corresponding conjectures. One of them is given below, see \cite{MT1973}, \cite{OD1987}.\\

\subsection{Pair Correlation Conjecture} For any real numbers $a, b$ with $0 \leq a < b < \infty$,
\begin{eqnarray}
	N(a,b,T)&=& \# \left \{ n\ne m:0 <\gamma_m,\gamma_n \leq T \text{ and } \frac{2 \pi a}{\log T} \leq \gamma_n-\gamma_m \leq \frac{2 \pi b}{\log T} \right \} \nonumber \\
	&\asymp & N(T) \int_a^b( 1-\sin(\pi x)^2) d x
\end{eqnarray}
as $T \to \infty$ .\\

\subsection{Zero Separation Conjecture} For all integers $n\leq 1$ and a small real number $\theta> 0$, the zeros spacing satisfy the
\begin{equation}
	\gamma_{n+1}-\gamma_n= O\left (\gamma_n^{-\theta} \right).
\end{equation}

This conjecture implies that the zeros of $\zeta(s)$ are simple, confer \cite{BB2006} for details and an application.\\

\subsection{Normalized Prime Gap Conjecture} The normalized prime gap has a Poisson distribution, specifically,
\begin{equation}
	\# \left \{ n \leq x:0\leq \frac{p_{n+1}-p_n}{\log p_n} \leq t \right \} =\left (1-e^{-t} \right ) x
\end{equation}
as $x \to \infty$, see \cite[p.\ 266]{IK2004}.\\

The above postulate is for small prime gaps $p_{n+1}-p_n=O(\log p_n)$. Accordingly, it seems that the distribution
of normalized prime gaps is the same as the Poisson distribution of normalized zero spacings.

\section{Prime Numbers Results} \label{sec3}
This section provides some elementary results on prime numbers, prime gaps, and related concepts. Let $p_n$ denotes the $n$th prime $p_n$ in the sequence of prime numbers $2, 3, 5, 7, \ldots,$ and let 
\begin{equation}
	\pi(x)= \#\{ p \leq x : p \text{ prime }\}
\end{equation}
denotes the prime counting function, respectively. \\

\begin{thm} {\normalfont(DelaValle Poussin 1899)}  \label{thm3.1} Let $x > 1$. Then 
	\begin{equation}
		\pi(x)= \li(x)+O\left ( xe^{-c \sqrt{\log x}} \right ) 
	\end{equation}
	where $c > 0$ is an absolute constant.
\end{thm}
\
The logarithm integral is defined by $\li(x)=\int_2^x \frac{1}{\log t}dt$. The proof of the strongest version $\pi(x)= \li(x)+O\left ( xe^{-c (\log x)^{3/5}(\log\log x)^{-1/5} } \right )$ of the prime number theorem is given in \cite[p.\ 307]{IV2003}, and the constant $c = .2018$ is given in \cite{FD2002}. Furthermore, the Riemann hypothesis claims the sharpest version $\pi(x)= \li(x)+O\left ( x^{1/2}\log^2 x \right )$.\\

\subsection{Formula for the $n$th Prime}
A formula for the $n$th prime in term of its position in the sequence of prime is stated here.\\

\begin{cor}  \label{cor3.1} Let $n$ be a large number and let $p_n$ be the $n$th prime. Then the followings statements are satisfied.\\
	\normalfont{(i)} The $n$th prime satisfies the inequalities 
	\begin{equation}
		an \log n \leq p_n \leq  bn \log n
	\end{equation} for some constants $a, b > 0$.\\
	\normalfont{(ii)} The $n$th prime satisfies  the asymptotic 
	\begin{equation}
		p_n =n\log n \left (1+O\left ( \frac{\log \log n}{ \log n}\right )\right )
	\end{equation} as the integer $n \to \infty$.
\end{cor}

This result is readily derived from the prime number theorem. A derivation of the complete infinite series appears in \cite{AJ2012}. For an arithmetic progression $q_n \equiv a \bmod q$, with $\gcd(a, q)=1$, the $n$th prime $q_n$ 
in the progression has the asymptotic formula $q_n \sim \varphi(q)n \log n$, see \cite{WS1979} for a
proof.\\

The prime number theorem, Theorem (\ref{thm3.1}), implies that the \textit{n}th prime has an asymptotic expression of the form \(p_n=n \log  n+O( n \text{loglog} n)\). The Cipolla formula specifies the exact formula for the \textit{n}th prime: \(p_n=(n \log  n)f( n) .\) The function \(f:\mathbb{N}\longrightarrow \mathbb{R}\) has the power series expansion  
\begin{eqnarray} \label{el31}
	f(n)&=&1+\frac{\log\log n-1}{\log  n}+\frac{\log\log
		n-2}{\log ^2 n}-\frac{\log ^2\log  n-6\log\log n+11}{\log ^3 n}\nonumber\\
	&&\hskip 2.35 in +\sum _{s\geq 3} \frac{f_s(\log\log n)}{\log ^{s+1} n},
\end{eqnarray}
where
\begin{equation}
	f_s(x)=a_sx^s+a_{s+1}x^{s-1}+\cdots +a_1x+a_0\in \mathbb{Q}[x]
\end{equation}
is a polynomial, \(\log ^s n=(\log  n)^s\) is an abbreviated notation for
the \textit{s}th power of the logarithm or its iterated form respectively. For example, the first few are the followings.\\
\begin{enumerate}
\item $	f_0(x)=x-1,$\item $ f_1(x)=x-2,$\item 
	$f_2(x)=-(x^2-6x+11)/2.$
\end{enumerate}

The leading coefficient \(a_s=\left.(-1)^{s+1}\right/s\) of the \textit{s}th polynomial \(f_s(x)\) has alternating sign for \(s\geq 1\). Other details on the asymptotic representation of the \textit{n}th prime \(p_n\) are discussed in \cite[p. \ 27]{AJ2012}, \cite{CM1902}, and \cite{DP1999}.\\

The asymptotic formula $p_n =n \log n +o(n \log n)$ probably can not be used to estimate the prime gap
$d_n=p_{n+1}-p_n$. This obstruction seems to arise from the fact that in the difference of two consecutive primes
$p_{n+1}-p_n=\log n +o(n \log n )$, which has an error term larger than the main term.

\subsection{Primes and Prime Gaps Results}
The prime gap $d_n=p_{n+1}-p_n$ is the
difference of two consecutive prime numbers, and the maximum prime gap is defined by
\begin{equation}
	d_{max}=\max_{p_n\leq x} \{p_{n+1}-p_n \}.
\end{equation}

The average prime gap of the short sequence $p_1, p_2, \ldots, p_n \leq x$ is given by the
asymptotic formula
\begin{equation}\label{el777}
	\overline{d_n}=\frac{x}{\pi(x)}= \log x \left (1+O(e^{\sqrt{\log x}}) \right ).
\end{equation}

These are important and intensively investigated local
properties of the prime numbers, see \cite{RN1996} for a survey.\\

The average statistic $\overline{d_n} \sim \log x$ for $p_n \leq x$ serves as guide on what to expect for the extreme values of the prime gaps. Naturally, the
interval $[1, x]$ contains a pair of primes with the average prime gap $p_{n+1}-p_n \approx \log p_n $, and the maximum
gap satisfies $p_{n+1}-p_n \geq \overline{d_n}$.\\

\subsection{Large Prime Gaps}
The simplest construction of an infinite sequence of primes with potentially large prime gaps is based on the factorial $n!$ numbers. In the double indexed infinite sequence consecutive integers
\begin{equation}
	n! + 1, \quad n! + 2, \quad \ldots, \quad n! + n, \quad n! + n + 1, \quad \ldots, 
\end{equation}
the integer $n!+1$ is a potential prime, but the block of $n-2$ consecutive integers $n!+2, \ldots, n! + n$ are
composites, and again, $n! + n + 1$ is a potential prime. The \textit{factorial primes} up to $q_n \leq x$ in this sequence have a
maximum prime gap of size $q_{n+1}- q_n \geq \log x /\log \log x$.\\

An improved version of the previous infinite sequence of numbers uses the \textit{primorial} numbers $p_n\# = 2\cdot 3\cdot 5\cdots p_n$. In this infinite sequence of integers
\begin{equation}
	p_n\# + 1, \quad p_n\# + 2, \quad \ldots, p_n\# + p_n, \quad p_n\# + p_n + 1, \quad \ldots, 
\end{equation}
the integer $p_n\# + 1$ is a potential prime, but the inner block of $n-2$ consecutive integers $p_n\# + 2, \ldots, p_n\# + p_n$ are composites, and again, $p_n\# + p_n + 1$ is a potential prime. The primorial primes up to $q_n \leq x$ in this sequence have
a maximum prime gap of size $q_{n+1}-q_n \approx \log x$.\\

The densities of most sequences of prime numbers (algebraic and nonalgebraic) are unknown, but can be
estimated using conditional Gaussian densities. The heuristic densities and the expected numbers of factorial
and primorial primes are expounded in \cite{CL2002}.\\

Despite the appearance, the prime gaps of the factorial and primorial primes are below or at the average prime
gap, viz, $d_n \leq \overline{d_n}$. The next level of complexity in the construction of an infinite sequence of
consecutive prime numbers with large prime gaps above the average value, viz, $d_n \geq \overline{d_n} \sim \log x$, is the
Westzinthius method. This method claims that there is an infinite sequence of primes $q_1, q_2, q_3, \ldots, q_n, \ldots,$ with gaps significantly larger than the average.\\

The most important technique for constructing an infinite sequence of consecutive prime numbers with large
prime gaps above the average value, viz, $d_n \geq  \log x$, is due to Westzynthius, see \cite{MW2007}.\\

\begin{thm} {\normalfont (Westzynthius)}  \label{thm3.2} Let $x > 1$ be a large real number, and let $p_n \leq x$ be the $n$th prime. Then
	\begin{equation}
		p_{n+1}-p_n\geq c \frac{(\log p_n) (\log \log \log p_n)}{\log \log \log \log p_n} , 
	\end{equation}
	where $c>0$ is a constant, for infinitely many $n\geq 1$ as $x \to \infty$. In particular, the limit supremum
	\begin{equation}
		\limsup_{n \to \infty} \frac{p_{n+1}-p_n}{ \log p_n}= \infty  
	\end{equation}
	holds. 
\end{thm}

There are several ways of proving this result, and all these techniques use both elementary methods and
advanced complicated methods, see \cite[p.\ 221]{MV2007}. The Westzynthius method has been developed by several
authors. These authors have reached concrete formulas and resolved the implied constant, see \cite{ER1935}, \cite{RA1938},
\cite{MP1990}, and \cite{PJ1997}.\\

\begin{thm} {\normalfont (\cite{PJ1997})} \label{thm3.3}  There is an infinite sequence of pairs of consecutive primes such that
	\begin{equation}
		p_{n+1}-p_n \geq c\frac{(\log p_n)(\log \log p_n) (\log \log \log \log p_n)}{(\log \log \log p_n)^2}  \nonumber
	\end{equation}
	for $p_n \leq x$, and $c = 2(e + o(1))$ is a constant.
\end{thm}

A probabilistic analysis of the prime numbers leads to the expression
\begin{equation}
	\limsup_{n \to \infty} \frac{p_{n+1}-p_n}{ \log^2 p_n}=1,
\end{equation}
which is known as Cramer conjecture. The current version of this conjecture calls for a larger limit of the form
$M = 2e^{-\gamma} > 1$, see \cite{GR1993}.\\

\subsection{Small Prime Gaps} 
The small prime gaps problem appears to be a much more difficult problem than the large prime gaps problem. The greatest lower bound of the sequence of small 
prime gaps coincides with the twin primes conjecture. The twin primes conjecture claims the existence of infinitely many pairs of primes such that 
$p_{n+1}- p_n = 2$.\\

\begin{thm} {\normalfont(\cite{GPY2007})}  \label{thm3.4}  The difference of consecutive primes satisfies 
	\begin{equation}
		\liminf_{n \to \infty} \frac{p_{n+1}-p_n}{ \sqrt{\log p_n}(\log \log p_n)^2}< \infty \nonumber .
	\end{equation}
\end{thm}

This is one of the closest approximation to the twin prime conjecture, it confirms the existence of infinitely
many primes pairs such that 
\begin{equation}
	p_{n+1}-p_n=O \left (\sqrt{\log p_n}(\log \log p_n)^2 \right ) .
\end{equation}
Clearly, this estimate is significantly smaller than the average prime gap, confer (\ref{el777}).\\

\section{Zeta Zeros Results} \label{sec4}
This section provides some elementary results on zeros, and zeros spacings of the zeta function as used in the proofs of the main results.
The most common data on the zeros of the zeta function is their density on the finite rectangle
\begin{equation}
	R(T) = \{ s \in \mathbb{C} : -T \leq \text{Im}(s) \leq T, \text{ and } -1 \leq \text{R} e(s) \leq 1 \}.
\end{equation}

The relevant formula for counting the number $N(T)$ of zeros of the zeta function on the rectangle in question
was proposed by Riemann and proved by vonMangoldt.\\

\begin{thm} {\normalfont (vonMangoldt)}  \label{thm4.1} Let $N(T) = \#\{ s \in R(T) : \zeta(s) = 0 \}$ be the number of complex zeros on the
	critical strip $0<\Re e(s) < 1$ such that $-T\leq \Im m(s) \leq T$. Then
	\begin{equation}
		N(T) = \frac{T}{2\pi }\log \frac{T}{2\pi } -\frac{T}{2\pi }+O(\log T)+S(T)\nonumber ,
	\end{equation}
	where $S(T)=\pi^{-1} \text{arg}( \zeta(1/ 2+ it))$ , and $T\geq 1$.
\end{thm}

\begin{proof}[\textbf{Proof}] The winding number enumerates the number of zeros of a function in a domain in the complex plane. Hence
	\begin{equation}
		N(T) =\frac{1}{2\pi } \int_{\delta R(T)}\frac{\xi^{'}(s)}{\xi(s) }ds
	\end{equation}
	enumerates the number of zeros of the entire function $\xi(s)=s(s-1)\pi^{-s/2}\Gamma(s/2)\zeta(s)$ for $s = \sigma+ it$ such that $|t| \geq 
	0 $ inside the rectangle $R(T)$. The integral easily splits into three line integrals:
	\begin{eqnarray}
		\frac{1}{i2\pi } \int_{\delta R(T)}\frac{\xi^{'}(s)}{\xi(s) }ds &=&  \frac{1}{i2\pi }  \left ( \int_{\delta R_1}\frac{\xi^{'}(s)}{\xi(s) }ds+ \int_{\delta R_2}\frac{\xi^{'}(s)}{\xi(s) }ds+ \int_{\delta R_3}\frac{\xi^{'}(s)}{\xi(s) }ds \right ) \nonumber \\
		&=&\frac{1}{2\pi }  \left ( \text{arg}_{\delta R_1} \xi(s) +{\normalfont arg}_{\delta R_2} \xi(s) +\text{arg}_{\delta R_3} \xi(s) \right ), 
	\end{eqnarray}
	where $\delta R_1 = [ 1 \to 2]$, $\delta R_2 =[2 \to 2 + iT \to  1/2 + iT]$, and $\delta R_3 =  [1/2 + iT \to 1 + iT \to 1]$ are the directed
	edges of the rectangle $R(T)$ or equivalently the paths of the line integrals, and the expression $\text{arg}_{\delta R_3} \xi(s)$ enumerates the changes in the argument of the function $\xi(s)$ along the curve $\delta R$. Since there is no change in the
	argument along the path, that is, $\text{arg}_{\delta R_1} \xi(s) =0$, the first line integral vanishes. The second and the third are
	equal since $\xi(s)=\xi (1- s)$, so 
	\begin{equation}
		N(T) = \frac{1}{2 \pi}  \text{arg}_{ \delta R_3 }  \xi (s)=   \frac{1}{ \pi} \left ( \text{arg} \pi^{-s/2} +  \text{arg} \Gamma (s/2)+ \text{arg} \zeta(1/2+it) \right )
	\end{equation}
	For more details see \cite[p.\ 19]{IV2003}, \cite[p.\ 127]{EH1974}, \cite[p.\ 160]{EL1985}.     
\end{proof}

The real valued function $S(t)=\pi^{-1} \text{arg}( \zeta(1/ 2+it) )$ along the path $\delta R_2$ has a prominent place in the analysis of zeta zeros. It is known that 
\begin{equation}
	S(t)= O(\log t) \quad \text{ and } \quad S(t ) /(2 \pi \sqrt{ \log  \log t} )
\end{equation}
has a normal distribution, see the literature.\\

\subsection{Formula for the $n$th Zero}
The zero counting theorem, Theorem \ref{thm4.1}, immediately leads to an asymptotic formula for the imaginary part of the $n$th zero.

\begin{cor}  \label{cor4.1}  The critical zeros $\rho_n =\beta_n+ i \gamma_n $ of $\zeta(s)$ satisfy the followings.\\
	\normalfont{(i)} The $n$th zero satisfies the inequalities 
	\begin{equation}
		\frac{an}{\log n }\leq \gamma_n \leq \frac{bn}{ \log n}\nonumber 
	\end{equation}
	for some constants $a, b > 0$.\\
	(ii) The imaginary part satisfies the asymptotic
	\begin{equation}
		\gamma_n = \frac{2 \pi n}{ \log n}\left (1+O\left(\frac{\log \log n}{ \log n} \right ) \right ) 
	\end{equation}
	as the integer $n \to \infty$.
\end{cor}

\begin{proof}[\textbf{Proof}] (ii) The inequality $N(\gamma_n- 1)<n \leq  N(\gamma_n)$ and the formula for $N(T)$ in Theorem \ref{thm4.1} imply that
	\begin{equation}
		n = \frac{\gamma_n}{2\pi }\log\left ( \frac{\gamma_n}{2 \pi} \right ) \left (1-\frac{2 \pi}{\log (\gamma_n/2 \pi) } +O(\frac{1}{\log \gamma_n})  \right )
	\end{equation}
	and
	\begin{equation}
		\log n =\log \left ( \frac{ \gamma_n}{2\pi } \right )\left (1+O\left ( \frac{\log \log \gamma_n}{\log \gamma_n} \right ) \right ).
	\end{equation}
	Combine these two items to confirm the claim. 
\end{proof}

Similar details on this result appear in \cite[p.\ 160]{EL1985}, \cite[p. 20]{IV2003}, et cetera. In \cite[p.\ 184]{NW2000} there is a proof for the
$n$th zero $\rho_n = \beta_n + i \gamma_n$ and its corresponding inequalities $an /\log n \leq | \rho_n| \leq bn / \log n$. \\

The zero counting theorem, Theorem (\ref{thm4.1}), implies that the \textit{n}th imaginary part of a zero has an asymptotic expression of the form $\gamma_n=\frac{2 \pi n}{ \log  n}+O( n \log \log n)$.
The expected exact formula for the \textit{n}th imaginary part should be of the form \(\gamma_n=(2 \pi n / \log  n)g( n),\) where the function \(g:\mathbb{N}\longrightarrow \mathbb{R}\) has a power series expansion \\ 
\begin{equation} \label{el35}
	g(n)=1+\sum _{s\geq 1} \frac{g_s(\log\log n)}{\log ^{s+1} n},
\end{equation}
and
\begin{equation}
	g_s(x)=b_sx^s+b_{s+1}x^{s-1}+\cdots +b_1x+b_0\in \mathbb{Q}[x]
\end{equation}\\
is a polynomial. The leading coefficient \(b_s=\left.(-1)^{s+1}\right/s\) of the \textit{s}th polynomial \(g_s(x)\) is expected to have alternating sign for \(s\geq 1\). \\

\subsection{Upper Bounds For Zero Gaps}  The difference $\gamma_{n+k}-\gamma_n=k/ \log(n+k)+o(n/ \log n)$, since the error term seems
to be larger than the main term. Thus, this result is probably not an effective way of estimating the size of the zero gap.

\begin{thm}   \label{thm4.2} {\normalfont (\cite{GG2007})} For all sufficiently large $n\geq 1$, the difference of every consecutive pair of zeros $\rho_{n}=\beta_{n}+i\gamma_{n}$ and $\rho_{n+1}=\beta_{n+1}+i\gamma_{n+1}$ of the
	zeta function $\zeta(s)$ on the critical line $\Re e(s)=1/2$ satisfie the inequality
	\begin{equation}
		\gamma_{n+1}-\gamma_n\leq \frac{\pi}{\log \log \gamma_n}(1+o(1))\leq c\nonumber ,
	\end{equation}
	where $c>0$ is a constant. \\
\end{thm}

\subsection{Simple Zero and Nonzero Gaps}  It is known that about $66\%$ of the zeros are simple.
More precisely, the number $N_1(T)=\#\{ 0\leq \mathcal{R} e(\rho)<1:\zeta(\rho)=0 \text{ and }\zeta^{`}(\rho)\ne0 \}$ of simple zeros $\rho= 1/2 + i\gamma$ of height $\gamma \leq T$ in the critical strip satisfies the
inequality $N_1(T) >(2/3 +o(1))(T / 2 \pi) \log T$, see \cite{MT1973}, \cite[p.\ 568]{IC2006} and similar papers.\\

\begin{cor}   \label{cor2.2} For a large integer $n \geq 1 $, the nth nontrivial zero of the zeta function is a simple zero.
\end{cor}

\begin{proof}[\textbf{Proof}] It is sufficient to show that the difference of every consecutive pair of zeros $\rho_{n}=\beta_{n}+i\gamma_{n}$ and $\rho_{n+1}=\beta_{n+1}+i\gamma_{n+1}$ has nonzero imaginary difference $\gamma_{n+1}-\gamma_n>0$. Toward this end, consider the strictly positive quantity $p_{n+1}-p_n\geq 2$. By Lemma \ref{lem5.1}, there is the lower bound
	\begin{eqnarray} 
		2 &\leq & p_{n+1}- p_n  \\
		&=&\frac{ \log ^2 n}{2 \pi} \left (1+O\left ( \frac{\log \log n}{ \log n}\right )\right ) \left (   \gamma_{n+k}-\gamma_n  \right ) \nonumber\\
		&\ll & \left ( \gamma_{n+1}-\gamma_n \right ) \log ^2 n  \nonumber. 
	\end{eqnarray} 
	Proceed to rearrange it as
	\begin{equation} 
		\frac{4 \pi}{\log ^2 n} \ll \left ( \gamma_{n+1}-\gamma_n   \right ) . 
	\end{equation} 
	
	Since $4 \pi /\log^2n >0$ for all integers $n\geq1$, (the smallest imaginary part is $\gamma_1 = 14.13472 \ldots$), this implies that all the zeros of the zeta function $\zeta(s)$ on the critical line $\Re e(s) = 1/2$ are simple. Otherwise, $4 \pi /\log^2n =0$ infinitely often, which is not plausible. 
\end{proof}


\section{Prime-Zero Duality Principle} \label{sec5}
The zero counting function (the Riemann-vonMangoldt formula) and the prime counting function (the prime
number theorem) are somewhat dual mathematical concepts. For large $T \geq 1$, these functions are given by
\begin{equation}
	N(T)=\frac{T}{ 2 \pi}\log  \left ( \frac{T}{ 2 \pi} \right )-\frac{T}{ 2 \pi} +O(\log T)
\end{equation}
and
\begin{equation}
	\pi(T)=\frac{T}{\log T}+O\left (\frac{T}{\log^2 T} \right),
\end{equation}
respectively. One of the striking similarities between these function are
the proofs, which are almost identical, and one of the acute differences between these function is the relatively
small error term $O(\log T)$ for the zeros counting function $N(T)$ compared to relatively large error term $O(T/\log^2T)$ for the primes counting function $\pi (T)$. But these
differences and similarities might have important applications. To understand the underpinning of this principle,
observe that Corollaries \ref{cor3.1} and \ref{cor4.1} lead to a simple spectral result, (this idea can be generalized to other pairs of
functions too). More precisely, a transformation from the primes domain to the zeros domain. Hence, certain
information on the zeros of the zeta function (more generally $L$-functions) has dual corresponding information
on the prime numbers and vice versa. Accordingly, the duality principle can be used to determine certain
unknown information on the primes from certain known information on the zeros and conversely.\\

A few of the dual properties of interest are described in a short prime-zero dictionary.
\begin{enumerate}
	\setlength{\itemindent}{-.15 in}
	\item The prime gap $p_{n+ 1}- p_n=\delta_1$ of two consecutive primes corresponds to the zero spacing $\gamma_{n+ 1}- \gamma_n=\epsilon_1$ of two consecutive zeros.
	\item  The sum of $k$  consecutive prime gaps $p_{n+ k}- p_n=(p_{n+ 1}- p_n)+\cdots +(p_{n+ k}- p_{n+k-1})=\delta_k $  corresponds to the sum of $v$ consecutive zero spacings $\gamma_{n+ k}- \gamma_n=(\gamma_{n+ 1}- \gamma_n)+\cdots +(\gamma_{n+ k}- \gamma_{n+k-1})=\epsilon_k$
	\item The prime gap $q_{n+ 1}- q_n=\delta_1$ of two consecutive primes $q_n, q_{n+1} \equiv a \bmod q$ in an arithmetic progression $\{ qm + a : \gcd(a, q) = 1, m \geq  0 \}$ corresponds to the zero spacing $\gamma_{n+ k}- \gamma_n=\epsilon_k$ of two consecutive zeros of the corresponding $L$-function.
\end{enumerate}
However, not every aspect of the functions $N(T)$ and $\pi(T)$ are well matched. Exempli gratia, a striking difference between $N(T)$ and $\pi(T)$ is the error terms. The former has a relatively small error term of logarithm size, but the latter has a large error term of exponential size.\\

\begin{lem} {\normalfont (Prime-Zero-Duality) } \label{lem5.1} For a sufficiently large integer $n \geq  1$, let $p_n$ be the $n$th prime and let $\rho_n = \beta_n
	+ i \gamma_n$ be the $n$th zero of the zeta function. Then, the followings hold independently of the small fixed integer $k < n$.\\
	\normalfont{(i)} The $n$th prime has the asymptotic formula
	\begin{equation}
		p_n= \frac{ \log ^2 n}{2 \pi}\left (1+O\left ( \frac{\log \log n}{ \log n}\right )\right ) \gamma_n.
	\end{equation}
	\normalfont{(ii)} The $n$th gap has the the asymptotic formula
	\begin{equation} 
		p_{n+k}- p_n= \left (   \gamma_{n+k}-\gamma_n  \right )\frac{ \log ^2 n}{2 \pi} \left (1+O\left ( \frac{\log \log n}{ \log n}\right )\right ) .
	\end{equation}
\end{lem}

\begin{proof}[\textbf{Proof}] (i) Let $k\geq 1$ be a small fixed integer, $1 \leq k < n$. By Corollaries \ref{cor3.1} and \ref{cor4.1}, the ratio $p_n/\gamma_n$ has the asymptotic formula
	\begin{eqnarray} \label{500}
		\frac{p_n}{\gamma_n}=\frac{\displaystyle n \log n \left (1+O\left ( \frac{\log \log n}{ \log n}\right )\right ) }{\displaystyle \frac{2 \pi n}{  \log n}\left (1+O\left ( \frac{\log \log n}{\log n} \right ) \right )}=\frac{ \log ^2 n}{2 \pi}\left (1+O\left ( \frac{\log \log n}{ \log n}\right ) \right ).
	\end{eqnarray}
	This prove the first claim. \\
	
	(ii) For a small fixed $k \geq 1$, the first term has this asymptotic:
	\begin{eqnarray} \label{505}
		p_{n+k}&=& \frac{ \log^2 (n+k)}{2 \pi} \left (1+O\left ( \frac{\log \log (n+k)}{ \log (n+k)}\right )\right ) \gamma_{n+k} \\
		&=& \frac{\log n}{2 \pi} \left (1+O \left (\frac{k}{n\log n} \right )  \right ) \left (1+O\left ( \frac{\log \log n}{ \log n}\right )\right ) \gamma_{n+k}\nonumber\\
		&=& \frac{\log^2 n}{2 \pi}  \left (1+O\left ( \frac{\log \log n}{ \log n}\right )\right )\gamma_{n+k} \nonumber ,
	\end{eqnarray}
	where
	\begin{equation}
		\log^2 (n+k)= (\log^2 n) \left (1+O \left (\frac{k}{n\log n} \right ) \right ) .
	\end{equation}
	Now, use (\ref{500}) and (\ref{500}) to compute the difference
	\begin{eqnarray}
		p_{n+k}-p_n&=& \frac{\log^2 (n+k)}{2 \pi}  \left (1+O\left ( \frac{\log \log (n+k)}{ \log (n+k)}\right )\right )\gamma_{n+k} \\
		&& \qquad  \qquad \qquad - \frac{\gamma_n}{2 \pi} \log^2 n\left (1+O\left ( \frac{\log \log n}{ \log n}\right )\right )\nonumber \\
		&=& \frac{\log^2 n}{2 \pi} \left (1+O\left ( \frac{\log \log n}{ \log n}\right )\right ) \left ( \gamma_{n+k}-\gamma_n\right ) \nonumber .
	\end{eqnarray}
\end{proof}

\begin{lem} {\normalfont (Prime-Zero-Duality)}  \label{lem5.2} For a sufficiently large integer $n \geq  1$, let $p_n$ be the $n$th prime and let $\rho_n = \beta_n
	+ i \gamma_n$ be the $n$th zero of the zeta function. Then, the followings hold independently of the small fixed integer $k < n$.\\
	\normalfont{(i)} The $n$th zero has the asymptotic formula
	\begin{equation}
		\gamma_n= \frac{2 \pi } {\log^2 n} \left (1+O\left ( \frac{\log \log n}{ \log n}\right )\right ) p_n.
	\end{equation}
	\normalfont{(ii)} The $n$th zero gap has the asymptotic formula
	\begin{equation}
		\gamma_{n+k}- \gamma_n=(p_{n+k}- p_n) \frac{4 \pi }{ \log ^2 n}\left (1+O\left ( \frac{\log \log n}{ \log n}\right )\right ) .
	\end{equation}
	
\end{lem}

\begin{proof}[\textbf{Proof}] (i) The ratio $\gamma_n / p_n $ has the asymptotic formula: 
	\begin{eqnarray}\label{550}
		\frac{\gamma_n}{p_n}=\frac{\displaystyle \frac{2 \pi n}{  \log n}\left (1+O\left ( \frac{\log \log n}{\log n} \right ) \right )}{\displaystyle n \log n \left (1+O\left ( \frac{\log \log n}{ \log n}\right )\right )       }=\frac{2 \pi} {\log^2 n} \left (1+O\left ( \frac{\log \log n}{ \log n}\right )\right ) .
	\end{eqnarray}
	This proves the first statement. (ii) For a small fixed $k \geq 1$, the first term has this asymptotic:
	\begin{eqnarray} \label{555}
		\gamma_{n+k}&=& \frac{2 \pi}{ \log^2 (n+k)} \left (1+O\left ( \frac{\log \log (n+k)}{ \log (n+k)}\right )\right )p_{n+k}  \\
		&=& \frac{2 \pi}{ \log^2 n}\left (1+O \left (\frac{k}{n\log n} \right )  \right ) \left (1+O\left ( \frac{\log \log n}{ \log n}\right )\right )p_{n+k} \nonumber \\
		&=& \frac{2 \pi}{ \log^2 n} \left (1+O\left ( \frac{\log \log n}{ \log n}\right )\right ) p_{n+k} \nonumber ,
	\end{eqnarray}
	
	where
	
	\begin{equation}
		\frac{1}{ \log^2 (n+k)}=\frac{1}{ \log^2 n}\left (1+O \left (\frac{k}{n\log n} \right )  \right ).
	\end{equation}
	
	Now, use (\ref{550}) and (\ref{555}) compute the difference
	\begin{eqnarray}
		\gamma_{n+k}- \gamma_n&=& \frac{2 \pi} {\log^2 (n+k)}  \left ( 1+O\left ( \frac{\log \log (n+k)}{ \log (n+k)} \right ) \right ) p_{n+k} \\
		&& \qquad  \qquad \qquad -\frac{2 \pi} {\log^2 n} \left ( 1+O\left ( \frac{\log \log n}{ \log n} \right ) \right ) p_{n}\nonumber \\
		&=&\left ( p_{n+k}-p_n \right )\frac{2 \pi}{ \log^2 n} \left (1+O\left ( \frac{\log \log n}{ \log n}\right )\right )  \nonumber ,
	\end{eqnarray}
	as claimed.
\end{proof}

\section{The Main Result} \label{sec6}

The notation $\log_k x= \log ( \cdots \log x \cdots )$ denotes the $k$-fold iterate of the logarithm function.
\begin{proof}[{ \textbf{Proof of Theorem \normalfont\ref{thm1.1}}}] (i) Let $n\geq 1$ be a large integer, and consider the strictly positive quantity $p_{n+1}-p_n\geq 2$. By Lemma \ref{lem5.1}, there is the lower bound
	\begin{eqnarray} 
		2&\leq &p_{n+k}- p_n= \left (   \gamma_{n+k}-\gamma_n  \right )\frac{ \log ^2 n}{2 \pi} \left (1+O\left ( \frac{\log \log n}{ \log n}\right )\right ) .
	\end{eqnarray}
	This implies that the right side is nonnegative. Rearranging it, and applying Theorem \ref{thm3.1} yield
	\begin{eqnarray} 
 c_0\frac{(\log p_n)(\log_2 p_n) (\log_4 p_n)}{(\log_3 p_n)^2}		&\leq & p_{n+1}- p_n  \\
	&=& 	\left (   \gamma_{n+k}-\gamma_n  \right )\frac{ \log ^2 n}{2 \pi} \left (1+O\left ( \frac{\log_2 n}{ \log n}\right )\right ) \nonumber.
	\end{eqnarray}
Clearly, the left side tends to infinity infinitely often as $n \to \infty$, with $c_0>0$ constant. Replacing $p_n \geq  c_1 n \log n \geq n$, see Lemma \ref{cor3.1}, and scaling both sides by $\log n$ yields
	\begin{eqnarray} 
c_0\frac{(\log_2 n) (\log_4 n)}{(\log_3 n)^2}		&\leq & p_{n+1}- p_n  \\
	&=& 	\left (   \gamma_{n+k}-\gamma_n  \right )\frac{ \log n}{2 \pi} \left (1+O\left ( \frac{\log_2 n}{ \log n}\right )\right )\nonumber
	\end{eqnarray}
	where $c_0, c_1, c_2, c_3, c_4>0$ are constants. Since the left side in the last inequality tends to infinity as $n \to \infty$, this proofs the limit supremum part in Theorem \ref{thm1.1}.\\
	
	(ii) Let $n\geq 1$ be a large integer, and consider the strictly positive quantity $p_{n+1}-p_n\geq 2$. By Lemma \ref{lem5.1}-ii, there is the lower bound
	\begin{eqnarray} 
		2&\leq & p_{n+1}- p_n \\
		&=& \left (   \gamma_{n+k}-\gamma_n  \right ) \frac{ \log^2 n}{2 \pi} \left (1+O\left ( \frac{\log \log n}{ \log n}\right )\right ) \nonumber.
	\end{eqnarray}
	This proves that the right side is nonnegative. Rearranging it, and applying Theorem \ref{thm3.2} yield
	\begin{eqnarray} 
		2 &\leq & \left (   \gamma_{n+k}-\gamma_n  \right )\frac{ \log^2 n}{2 \pi} \left (1+O\left ( \frac{\log_2 n}{ \log n}\right )\right )  \\
		&=&p_{n+1}- p_n  \nonumber\\
		&\leq & c_3\sqrt{\log p_n}(\log_2 p_n)^2 \nonumber
	\end{eqnarray}
	infinitely often as $n \to \infty$, with $c_3>0$ constant. Replacing $p_n \leq  2n \log n$, see Lemma \ref{cor3.1}, and scaling both sides by $\log n$ yield
	\begin{eqnarray}  
		\frac{2}{\log n} &\leq & \left (   \gamma_{n+k}-\gamma_n  \right ) \frac{ \log n}{2 \pi} \left (1+O\left ( \frac{\log_2 n}{ \log n}\right )\right )  \\
		&\leq & c_3 \frac{(\log_2 p_n)^2}{\sqrt{\log p_n}}  \nonumber\\
		&\leq & c_4 \frac{(\log_2 n)^2}{\sqrt{\log n}}\nonumber.
	\end{eqnarray}
	Since the right side in the last inequality tends to zero as $n \to \infty$, this proofs the limit infimum part in Theorem \ref{thm1.1}.
\end{proof}



\currfilename.\\

\end{document}